\pdfoutput=1
\RequirePackage{ifpdf}
\ifpdf 
\documentclass[pdftex]{sigma}
\else
\documentclass{sigma}
\fi

{ \theoremstyle{definition}
\newtheorem{Example}{Example}[section]
}
\numberwithin{equation}{section}

\begin{document}

\allowdisplaybreaks

\newcommand{\arXivNumber}{1709.09682}

\renewcommand{\thefootnote}{}

\renewcommand{\PaperNumber}{003}

\FirstPageHeading

\ShortArticleName{Manifold Ways to Darboux--Halphen System}

\ArticleName{Manifold Ways to Darboux--Halphen System\footnote{This paper is a~contribution to the Special Issue on Modular Forms and String Theory in honor of Noriko Yui. The full collection is available at \href{http://www.emis.de/journals/SIGMA/modular-forms.html}{http://www.emis.de/journals/SIGMA/modular-forms.html}}}

\Author{John Alexander Cruz MORALES~$^{\dag^1}$, Hossein MOVASATI~$^{\dag^2}$, Younes NIKDELAN~$^{\dag^3}$,\\ Raju ROYCHOWDHURY~$^{\dag^4}$ and Marcus A.C.~TORRES~$^{\dag^2}$}

\AuthorNameForHeading{J.A.C.~Morales, H.~Movasati, Y.~Nikdelan, R.~Roychowdhury and M.A.C.~Torres}

\Address{$^{\dag^1}$~Departamento de Matem\'aticas, Universidad Nacional de Colombia, Bogot\'a, Colombia}
\EmailDD{\href{mailto:jacruzmo@unal.edu.co}{jacruzmo@unal.edu.co}}

\Address{$^{\dag^2}$~Instituto Nacional de Matem\'atica Pura e Aplicada (IMPA), Rio de Janeiro, Brazil}
\EmailDD{\href{mailto:hossein@impa.br}{hossein@impa.br}, \href{mailto:mtorres@impa.br}{mtorres@impa.br}}
\URLaddressDD{\url{http://w3.impa.br/~hossein/}}

\Address{$^{\dag^3}$~Instituto de Matem\'atica e Estat\'{i}stica (IME),\\
\hphantom{$^{\dag^3}$}~Universidade do Estado do Rio de Janeiro (UERJ), Rio de Janeiro, Brazil}
\EmailDD{\href{mailto:younes.nikdelan@ime.uerj.br}{younes.nikdelan@ime.uerj.br}}
\URLaddressDD{\url{https://sites.google.com/site/younesnikdelan/}}

\Address{$^{\dag^4}$~Instituto de F\'isica, Universidade de S\~{a}o Paulo (IF-USP), S\~{a}o Paulo, Brazil}
\EmailDD{\href{mailto:raju@if.usp.br}{raju@if.usp.br}}

\ArticleDates{Received September 29, 2017, in f\/inal form January 03, 2018; Published online January 08, 2018}

\Abstract{Many distinct problems give birth to Darboux--Halphen system of dif\/ferential equations and here we review some of them. The f\/irst is the classical problem presented by Darboux and later solved by Halphen concerning f\/inding inf\/inite number of double orthogonal surfaces in~$\mathbb{R}^3$. The second is a problem in general relativity about gravitational instanton in Bianchi~IX metric space. The third problem stems from the new take on the moduli of enhanced elliptic curves called \textit{Gauss--Manin connection in disguise} developed by one of the authors and f\/inally in the last problem Darboux--Halphen system emerges from the associative algebra on the tangent space of a Frobenius manifold.}


\Classification{34M55; 53D45; 83C05}

\renewcommand{\thefootnote}{\arabic{footnote}}
\setcounter{footnote}{0}

\section{Introduction}\label{intro}
The Darboux--Halphen system of dif\/ferential equations
\begin{gather}
\dot t_1=t_1(t_2+t_3)-t_2t_3,\nonumber\\
\dot t_2= t_2(t_1+t_3)-t_1t_3, \nonumber \\
 \dot t_3= t_3(t_1+t_2)-t_1t_2, \qquad \dot {}=\partial /\partial \tau,\label{halphen}
\end{gather}
where $\tau$ is a free parameter, f\/irst came to existence when Darboux \cite{MR1508661} was studying the existence of an inf\/inite number of double orthogonal system of coordinates. He formulated the problem as follows: {\it Let~$A$ and~$B$ be two f\/ixed surfaces in the $3$-dimensional Euclidean space~$\mathbb{R}^3$ and suppose that~$\Sigma$ is the family of surfaces which are the locus of the points such that the sum of their distances from the surfaces~$A$ and~$B$ are constant; and~$\Sigma'$ is the family of surfaces which are the locus of the points so that the difference of their distances from the surfaces~$A$ and~$B$ are constant. Is there a third family of surfaces intersecting~$\Sigma$ and~$\Sigma'$ orthogonally?} When we restrict the third family to the surfaces given by second degree equations, we f\/ind the Darboux--Halphen system. In Section~\ref{darbouxproblem} we present Halphen's solution to this problem.

The Darboux--Halphen system also emerge from a direct map from Ramanujan relations (Section~\ref{Ramanujan relations}).

In 1979, Gibbons and Pope \cite{MR535151} found the Darboux--Halphen equations while studying gra\-vitational instanton solutions in Bianchi~IX spaces without having noticed it. Couple of decades later, Ablowitz et al.~\cite{MR1713583} pointed it out and recently one of the authors \cite{2016IJGMM..1350042C} explored its integrability aspects. A gravitational instanton is simply the (anti-)self-duality condition imposed on the curvature of a Einstein manifold with asymptotic locally Euclidian boundary conditions.

Hitchin \cite{hitchin1995twistor} and Tod~\cite{tod1994self} realized that (anti-)self-duality in Bianchi~IX metric has a more general solution envolving to a Darboux--Halphen system coupled to another system of linear dif\/ferential equation similar to Darboux--Halphen. A revised and simplif\/ied proof of the results of Tod and Hitchin can be found in \cite{babich1998self}. See \cite{manin2015symbolic} for a physical application in cosmology. We review these works in Section~\ref{BianchiIXsection}.

Another author \cite{ho06-1,ho14} among us met the Darboux--Halphen system while exploring the Gauss--Manin connection of a universal family of elliptic curves. This method is called \textit{Gauss--Manin connection in disguise}, which also name the vector f\/ield in this method that gives rise to the Darboux--Halphen equations and we present it in Section~\ref{GMCD}.

The last interesting problem where Darboux--Halphen system appears is in the context of a~3-dimensional Frobenius manifold with a certain potential function~$F(t)$. Frobenius manifold arose as a geometrization of Witten--Dijkgraaf--Verlinde--Verlinde (WDVV) equations \cite{Dijkgraaf:1990tp,Witten:1990fz}, an overdetermined system of dif\/ferential equations that appear in the physics of topological f\/ield theories in 2 dimensions. In this particular case of dimension 3, the WDVV equation is known as Chazy equation, which has a~close tie with the solutions of Darboux system. We present it in Section~\ref{Frobenius}, where we follow Dubrovin's notes~\cite{Dubrovin1992NuPhB,Dubrovin:1996di}.

We conclude this article, crossing information between problems displayed here, which led us to interesting remarks and an evidence that leads to a new way on how to examine spectral curves from monopoles using Gauss--Manin connection in disguise, further explored in~\cite{Torres:2017hqc} by one of the authors.

Throughout this text we make extensive use of Einstein summation convention where the sum over identical upper and lower indices is implicit.

\section{The Darboux problem}\label{darbouxproblem}
The above Darboux problem given in Section \ref{intro} is equivalent to the following problem: \emph{Let~$A$ and~$B$ be as before and suppose that $\Sigma$ is a family of surfaces parallel to~$A$ which is parameterized by~$v$, and~$\Sigma '$ is a family of surfaces parallel to $B$ that is parameterized by~$w$. Is there a third family of surfaces parameterized by $\tau$ such that it intersects $\Sigma$ and $\Sigma '$ orthogonally?} Note that, two surfaces $A_1$ and $A_2$ are said to be \emph{parallel}, if there exist a constant $c\in \mathbb{R}_{\ne 0}$ and a continuous one to one map between points $a_1\in A_1$ and points $a_2\in A_2$, such that the tangent planes at these points are parallel and the position vector $\bf a_2=a_1+{\rm c}\hat{N}$, where $\bf \hat{N}$ is the unitary vector normal to the surface $A_1$ at~$a_1$. We say that a~family of surfaces is parameterized by $\varphi=\varphi(x,y,z)$, if any surface belonging to this family is given by $\varphi(x,y,z)={\rm const}$, in which $x$, $y$, $z$ are the standard coordinates of~$\mathbb{R}^3$. If for a~function $\varphi=\varphi(x,y,z)$, we def\/ine
 \begin{gather*}
\varphi_x=\frac{\partial \varphi}{\partial x} , \qquad
\varphi_y=\frac{\partial \varphi}{\partial y} , \qquad
\varphi_z=\frac{\partial \varphi}{\partial z} ,
\end{gather*}
then in the latter problem, Darboux chose a case of parametrization of parallel surfaces that gives the Gauss map for points on the parallel surfaces
\begin{gather*}
v_x^2+v_y^2+v_z^2=1,\qquad w_x^2+w_y^2+w_z^2=1,
\end{gather*}
and the condition of orthogonality at points in the intersection of $v$ and $\tau$ and $w$ and $\tau$, respectively, is given by
\begin{gather*}
\tau_xv_x+\tau_yv_y+\tau_zv_z=0,\qquad \tau_xw_x+\tau_yw_y+\tau_zw_z=0.
\end{gather*}
So the problem is equivalent to the following system of equations,
\begin{gather}
v_x^2+v_y^2+v_z^2=1, \label{eq 1}\\
w_x^2+w_y^2+w_z^2=1, \label{eq 2}\\
\tau_xv_x+\tau_yv_y+\tau_zv_z=0, \label{eq 3}\\
\tau_xw_x+\tau_yw_y+\tau_zw_z=0. \label{eq 4}
\end{gather}
If for a function $\varphi$ of three variables $(x,y,z)$ we def\/ine the operator
\begin{gather*}\partial_\varphi:=\varphi_x\frac{\partial}{\partial
x}+\varphi_y\frac{\partial}{\partial y}+\varphi_z\frac{\partial}{\partial z},\end{gather*} then equations \eqref{eq
1}, \eqref{eq 2}, \eqref{eq 3} and \eqref{eq 4}, respectively, are given by $\partial_vv=1$, $\partial_ww=1$, $\partial_\tau v=\partial_v\tau =0$ and $\partial_\tau w=\partial_w\tau =0$, respectively. These equations imply $\partial_\tau \partial_vv=0$, $\partial_\tau \partial_ww=0$, $\partial_v\partial_v\tau =0$, $\partial_w\partial_w\tau =0$. Hence we get
\begin{gather} \label{eq parv}
2\partial_v\partial_v\tau -\partial_\tau \partial_vv=0,\qquad
2\partial_w\partial_w\tau -\partial_\tau \partial_ww=0.
\end{gather}

The situation is more interesting when the family $(\tau )$ is of second degree. Hence let us suppose that the family $(\tau )$ is
given by
\begin{gather}\label{eq elips}
ax^2+by^2+cz^2=1,
\end{gather}
where $a$, $b$, $c$ are functions of the parameter $\tau $. By this assumption, equations \eqref{eq 3} and \eqref{eq 4} yield
\begin{gather}\label{eq axv}
 axv_x+byv_y+czv_z=0,\qquad
 axw_x+byw_y+czw_z=0. 
\end{gather}
As well, from equation \eqref{eq parv} we get
\begin{gather}\label{eq avx}
av_x^2+bv_y^2+cv_z^2=0,\qquad
aw_x^2+bw_y^2+cw_z^2=0.
\end{gather}
Equations \eqref{eq 1}, \eqref{eq elips}, \eqref{eq axv} and \eqref{eq avx} imply
\begin{gather}
\big(a^2b'+b^2a'\big)(xv_y-yv_x)^2+\big(b^2c'+c^2b'\big)(yv_z-zv_y)^2 +\big(c^2a'+a^2c'\big)(zv_x-xv_z)^2=0, \nonumber\\
ab(xv_y-yv_x)^2+bc(yv_z-zv_y)^2+ca(zv_x-xv_z)^2=0,\label{eq v}
\end{gather}
in which $'=\frac{d}{d\tau }$. Analogously for $w$ we f\/ind
\begin{gather}
\big(a^2b'+b^2a'\big)(xw_y-yw_x)^2+\big(b^2c'+c^2b'\big)(yw_z-zw_y)^2 +\big(c^2a'+a^2c'\big)(zw_x-xw_z)^2=0,\nonumber\\
ab(xw_y-yw_x)^2+bc(yw_z-zw_y)^2+ca(zw_x-xw_z)^2=0. \label{eq w}
\end{gather}
The two equations in \eqref{eq v} and the two equations in \eqref{eq w} become equivalent if
\begin{gather}\label{eq dar1}
\frac{a^2b'+b^2a'}{ab}=\frac{b^2c'+c^2b'}{bc}=\frac{c^2a'+a^2c'}{ac}.
\end{gather}
If in \eqref{eq elips} we substitute $a$, $b$, $c$ respectively by $\frac{1}{t_1}$, $\frac{1}{t_2}$, $\frac{1}{t_3}$, then from~\eqref{eq
dar1} we get that the family
\begin{gather*}
\frac{x^2}{t_1}+\frac{y^2}{t_2}+\frac{z^2}{t_3}=1,
\end{gather*}
is orthogonal to both $\Sigma$ and $\Sigma'$ if $t_1$, $t_2$, $t_3$
satisfy the following
\begin{gather}\label{eq darboux}
t_3\left(\frac{{\rm d}t_1}{{\rm d}\tau }+\frac{{\rm d}t_2}{d\tau }\right)=t_2\left(\frac{{\rm d}t_1}{{\rm d}\tau }+\frac{{\rm d}t_3}{{\rm d}\tau }\right)=t_1\left(\frac{{\rm d}t_2}{{\rm d}\tau }+\frac{{\rm d}t_3}{{\rm d}\tau }\right).
\end{gather}
A particular case of the equation \eqref{eq darboux}, which is known
as Darboux--Halphen system, is given in \eqref{halphen}.
In 1881, G.~Halphen~\cite{ha81} studied this system of dif\/ferential
equations and expressed a~solution of it in terms of the logarithmic
derivatives of the theta functions; namely,
\begin{gather}
t_1= 2 (\ln \theta_2(\tau ))',\qquad
t_2=2 (\ln \theta_3(\tau ))', \qquad
t_3=2 (\ln \theta_4(\tau ))'.\label{darbouxsolution}
\end{gather}
with
\begin{gather*}
\theta_2(\tau ):=\sum_{n=-\infty}^\infty q^{\frac{1}{2}(n+\frac{1}{2})^2},\qquad
\theta_3(\tau ):=\sum_{n=-\infty}^\infty q^{\frac{1}{2}n^2},\\
\theta_4(\tau ):=\sum_{n=-\infty}^\infty (-1)^nq^{\frac{1}{2}n^2},\qquad q=e^{2\pi i \tau },\qquad \tau \in {\mathbb H}.
\end{gather*}
These theta functions can be written in terms of the more general theta functions with characteristics $r$ and $s$ and arguments $z$ and $\sigma$:
\begin{gather*}
\vartheta[r,s](z,\sigma)=\sum_{m\in \mathbb{Z}}\exp\big\{\pi i(m+r)^2\sigma+2\pi i (m+r)(z+s)\big\},\qquad z,r,s\in \mathbb{C},\qquad \sigma \in \mathbb{H},
\end{gather*}
such that
\begin{gather*}
\theta_2(\tau )= \vartheta[1/2,0](0,\tau ),\qquad\theta_3(\tau )= \vartheta[0,0](0,\tau ),\qquad\theta_4(\tau )= \vartheta[0,1/2](0,\tau ).
\end{gather*}

\section{Ramanujan relations between Eisenstein series}\label{Ramanujan relations}
The following dif\/ferential equation
\begin{gather}\label{raman}
q\frac{\partial E_2}{\partial q}=\frac{1}{12}\big(E_2^2-E_4\big), \qquad
q\frac{\partial E_4}{\partial q}=\frac{1}{3}(E_2E_4-E_6), \qquad
q\frac{\partial E_6}{\partial q}=\frac{1}{2}\big(E_2E_6-E_4^2\big),
\end{gather}
where $E_i$'s are the Eisenstein series
\begin{gather*}
E_{2i}(q):=1+b_i\sum_{n=1}^\infty \left (\sum_{d\mid n}d^{2i-1}\right )q^{n},\qquad i=1,2,3,
\end{gather*}
and $(b_1,b_2,b_3)=(-24, 240, -504)$, was discovered by Ramanujan in \cite{ra16} and it is mainly known
as Ramanujan's relations between Eisenstein series.
Ramanujan was a master of formal power series
and had a very limited access to the modern mathematics of his time.
In particular, he and many people in number theory didn't know that the dif\/ferential equation~\eqref{raman} had already been studied by Halphen in his book \cite[p.~331]{hal00}, thirty years before S.~Ramanujan.
The equalities of the coef\/f\/icients of $x^i$ in
\begin{gather}\label{27spet2017}
4(x-t_1)(x-t_2)(x-t_3)=4(x-a_1E_2)^3-a_2E_4(x-a_1E_2)-a_3E_6,
\end{gather}
where
\begin{gather*}
(a_1,a_2,a_3):=\left(\frac{2\pi i}{12}, 12\left(\frac{2\pi i}{12}\right)^2, 8\left(\frac{2\pi i}{12}\right)^3\right),
\end{gather*}
gives us a map from ${\mathbb C}^3$ into itself which transforms Darboux--Halphen into Ramanujan dif\/fe\-ren\-tial equation,
see \cite[pp.~330, 335]{ho14}.

\section{Self-duality in Bianchi IX metrics}\label{BianchiIXsection}
An instanton is a f\/ield conf\/iguration that vanishes at spacetime inf\/inity. It is the quantum ef\/fect that leads metastable states to decay into vacuum. It is a~phenomenon that takes place in usual spacetime with signature $(-,+,+,+)$ but in order to perform physical calculation we use its equivalence with a soliton solution (static and energetically stable f\/ield conf\/iguration) in Euclidean spacetime.
 In Yang--Mills theory, self-duality of the f\/ield strength $F_{\mu\nu}=\epsilon_{\mu\nu\rho\sigma}F^{\rho\sigma}$ in four spacetime dimensions is a widely known instanton conf\/iguration \cite{belavin1975pseudoparticle}. Similarly, self-duality constraint on the curvature two-form (and connection $1$-form) in Cartan's formalism of general relativity characterizes a gravitational instanton. An important feature of self-duality of the curvature is that the Ricci-tensor vanishes and it is a solution of the vacuum Einstein equations.
Also, self-dual curvature leads to solving a linear dif\/ferential equation, a task much easier than solving the full non-linear Einstein equations. Gravitational instantons were found in Bianchi~IX metrics, by Gibbons and Pope~\cite{MR535151}. Without realizing it, they arrived at Darboux--Halphen system from self-duality constraints.

In \cite{Bianchi2001}, L.~Bianchi studied continuous isometries of 3-dimensional spaces. He noticed that the continuous isometries (continuous motion that preserve ${\rm d}s^2$) of a space form a f\/inite-dimensional Lie group and he classif\/ied such spaces according to the corresponding group of isometries. Bianchi~IX corresponds to a~3-dimensional space with ${\rm SO}(3)$ or ${\rm SU}(2)$ as Lie group of isometries. When we consider it in the context of 4-dimensional cosmology, the isometries lie in the 3 spacial directions \cite{manin2015symbolic}, but since we are working in Euclidean signature we consider the isometry group ${\rm SO}(3)$ as a subgroup of ${\rm SO}(4)$. In this conf\/iguration, as the instanton vanishes at inf\/inity, Lorentz symmetry is recovered and the space is called asymptotically locally Euclidean (ALE). This same manifold describes the reduced\footnote{Moduli of charge 2 monopoles reduced by quotient by $\mathbb{R}^3$ action.} moduli~$M^0_2$ of charge 2 monopoles in a~${\rm SU}(2)$ Yang--Mills--Higgs theory.

A magnetic 2-monopole is a soliton solution of charge 2 of Bogomolny equations in the Yang--Mills--Higgs theory in $\mathbb{R}^3$, where ${\rm SU}(2)$ Yang--Mills is a gauge theory of 1-form connections $\bf{A}$ on a principal ${\rm SU}(2)$-bundle while the Higgs f\/ield $\Phi$ correspond to a section of an associated $\mathfrak{su}(2)$-bundle \cite{atiyahgeometry,Hurtubise:1983hl}. In \cite{atiyahgeometry}, Atiyah and Hitchin showed that the reduced moduli $M^0_2$ of 2-monopoles is a~4-dimensional hyperk\"ahler manifold and an anti-self-dual (curvature-wise) Einstein manifold. Since $M^0_2$ admits ${\rm SO}(3)$ isometry, the metric is a Bianch IX\footnote{Note that here the four coordinates of the moduli are not spacetime directions, but internal parameters of a~2-monopole solution.} (\ref{bianchiix}). This is a consequence of the hyperk\"ahler structure of $M^0_2$ which has an $S^2$-parameter family of complex structures, i.e., if $I$, $J$, $K$ are covariantly constant complex structures in $M^0_2$ then $aI+bJ+cK$ is also a covariantly
constant complex structure in $M^0_2$ given that $a^2+b^2+c^2=1$.

Here we present a detailed derivation of the Darboux--Halphen system starting from the Euclidean Bianchi~IX metric with ${\rm SO}(3)$ symmetry with an imposition of the constraints of self-duality at the level of Riemann curvature. The constraint of anti-self-duality yields an anti-instanton, a solution with negative instanton number and we present this solution together by using $\pm \operatorname{sign}$.
We follow the steps of \cite{2016IJGMM..1350042C} and \cite{MR535151}, see also \cite{PV2012}.

\subsection{Geometric analysis}
A metric for a 4-dimensional spacetime with coordinates $(x^1,x^2,x^3,x^4)$, Euclidean time~$x^4$
and ${\rm SO}(3)\subset {\rm SO}(4)$ isometry is written in terms of invariant 1-forms $\sigma^i$ on ${\rm SO}(3)$, dual to the standard basis $X_1$, $X_2$, $X_3$ of its Lie algebra
\begin{gather*}
\sigma^i = - \frac1{r^2} \eta^{i}_{\mu \nu} x^\mu {\rm d}x^\nu,
\end{gather*}
where $\mu, \nu = 1,2,3,4$ and $i= 1,2,3$ and $\eta^{i}_{\mu \nu}$ is a 't Hooft symbol given by
\begin{gather*}
\eta^i_{\mu \nu} = \begin{cases}
\varepsilon_{i\mu\nu}+\delta_{i\mu}\delta_{\nu 4}-\delta_{\mu 4}\delta_{i\nu}\qquad\text{or},\\
\varepsilon_{i\mu\nu}-\delta_{i\mu}\delta_{\nu 4}+\delta_{\mu 4}\delta_{i\nu}.
\end{cases}
\end{gather*}
according to two dif\/ferent choices of $\mathfrak{so}(3)$ generators in the Lie algebra of the group ${\rm SO}(4)= {\rm SU}(2)\times {\rm SU}(2)$. Among the symbols presented above, $\varepsilon_{i\mu\nu}$ is the Levi-Civita symbol
\begin{gather*}
\varepsilon_{123}=\varepsilon_{231}=\varepsilon_{312}=-\varepsilon_{213}=-\varepsilon_{321}=-\varepsilon_{132}=1,\qquad\text{and zero elsewhere,}
\end{gather*}
and $\delta_{i\mu}$ refers to the Kronecker delta. The $\sigma^i$'s obey the structure equation:
\begin{gather}\label{struc} {\rm d} \sigma^i = - {\varepsilon}_{ijk} \sigma^j \wedge \sigma^k , \end{gather}
where we use Einstein summation in the repeated upper and lower indices here and what follows below. This choice of ${\rm SO}(3)$ isometry leads to a 4D spherically symmetric Bianchi IX metric			
\begin{gather}
\label{bianchiix} {\rm d}s^2 = c_0(r)^2 {\rm d}r^2 + c_1^2 (r) \big(\sigma^1\big)^2 + c_2^2 (r) \big(\sigma^2\big)^2 + c_3^2 (r) \big(\sigma^3\big)^2 ,
\end{gather}
with $r=\sqrt[]{x_1^2+x_2^2+x_3^2+x_4^2}$, $ c_0 (r) = c_1 (r) c_2 (r) c_3 (r)$ and $c_1$, $c_2$, $c_3$ being functions of $r$.

We can impose self-duality in Bianchi IX metric in two ways:
\begin{enumerate}\itemsep=0pt
\item[1)] connection wise self-duality,
\item[2)] curvature wise self-duality.
\end{enumerate}
The connection wise self-duality is a stronger form of self-duality that leads to self-dual curvature tensor \cite{Eguchi:1978bx}. This form of self-duality does not present Darboux--Halphen system, but the Lagrange or Euler-top system~\cite{2016IJGMM..1350042C}. It is not in our goal to describe it here.

 One can perform a standard analysis using vierbeins, leading to Cartan's structure equation. The vierbeins could be chosen as
 \begin{gather*}
 e^0 = c_0 {\rm d}r, \qquad e^i = c_i \sigma^i \quad (\textnormal{no sum in}\ i), \qquad i = 1, 2, 3,
 \end{gather*}
 and the connection $1$-form can be obtained from the structure equation
 \begin{gather*}
{\rm d}e^a= e^b\wedge {\omega^a}_b,
 \end{gather*}
 where $a, b=0,1,2,3$. Obviously, $e^0$ produces no connections while other three does
 \begin{gather}
 {\rm d} e^0 = 0,\qquad
 {\rm d} e^i = \partial_r c_i {\rm d}r \wedge \sigma^i - c_i {\varepsilon}_{ijk} \sigma^j \wedge \sigma^k. \label{ei}
 \end{gather}
The f\/irst term on the r.h.s.\ above gives ${\omega^i}_0$ while the second term needs to be rewritten in order to produce a antisymmetric connection $1$-form
 \begin{gather*}
 \varepsilon_{ijk}\frac{c_i^2}{c_i} \sigma^j \wedge \sigma^k = \varepsilon_{ijk}\frac{2c_i^2 + \big( c_j^2 - c_k^2 \big) - \big( c_j^2 - c_k^2 \big)}{ 2c_i} \sigma^j \wedge \sigma^k \nonumber \\
\hphantom{\varepsilon_{ijk}\frac{c_i^2}{c_i} \sigma^j \wedge \sigma^k}{}
= \varepsilon_{ijk} \frac{c_i^2 + c_j^2 - c_k^2 }{ 2c_i c_j} e^j \wedge \sigma^k + \varepsilon_{ikj}\frac{c_i^2 + c_k^2 - c_j^2 }{ 2c_i c_k} e^k \wedge \sigma^j\\
\hphantom{\varepsilon_{ijk}\frac{c_i^2}{c_i} \sigma^j \wedge \sigma^k}{}
 = \varepsilon_{ijk} \frac{c_i^2 + c_j^2 - c_k^2 }{ c_i c_j} e^j \wedge \sigma^k .
 \end{gather*}
 Rewriting (\ref{ei}),
 \begin{gather*} 
{\rm d} e^i = - \frac{\partial_r c_i}{c_0} \sigma^i \wedge e^0 - \varepsilon_{ijk} \frac{c_i^2 + c_j^2 - c_k^2 }{ c_i c_j} e^j \wedge \sigma^k .
 \end{gather*}
Hence,
 \begin{gather}\label{connection1form}
 {\omega^i}_0 = \frac{\partial_r c_i}{c_0} \sigma^i \quad (\text{no sum in} \ i) ,\qquad {\omega^i}_j = - \varepsilon_{ijk} \frac{c_i^2 + c_j^2 - c_k^2 }{ c_i c_j} \sigma^k.
 \end{gather}
Here the connection 1-form components are anti-symmetric under permutation of its indices.

\subsection{Curvature wise self-duality and Darboux--Halphen system}

Curvature-wise self-duality was f\/irst studied in search of gravitational instantons. It is a more general solution than imposing self-duality on connection 1-forms. The Cartan-structure equation for Ricci tensor is
\begin{gather*}\label{cartan2}
 R_{ij} = {\rm d} \omega_{ij} + \omega_{im} \wedge {\omega^m}_j.
\end{gather*}
The (anti-)self-duality of curvature demands that
\begin{gather}
\label{curvsd} R_{0i} =\pm \frac12 { \varepsilon_{0ilm} } R^{lm} = \pm R_{jk},
\end{gather}
where $\{i,j,k\}$, in this order, are a cyclic permutation of $\{1,2,3\}$ and we used the fact that Euclidean vierbein indices are raised and lowered with Kronecker deltas ${\delta^i}_j$. Comparing the l.h.s.\ and r.h.s.\ of~(\ref{curvsd}), we have
\begin{gather}
{\rm d}(\omega_{0i}\mp\omega_{jk})=\pm(\omega_{0k}\mp\omega_{ij})\wedge(\omega_{0j}\mp\omega_{ki}),\nonumber\\
{\rm d}\big(\lambda_1(r)\sigma^i\big)=\pm\big(\lambda_3(r)\sigma^k\big)\wedge\big(\lambda_2(r)\sigma^j\big)=\lambda_3(r)\lambda_2(r)\sigma^k\wedge\sigma^j, \nonumber\\
\partial_r\lambda_1 dr\wedge\sigma^i+\lambda_1d\sigma^i=\mp\lambda_2\lambda_3\sigma^j\wedge\sigma^k,\label{sdR}
\end{gather}
where the second line comes from equation (\ref{connection1form}) with $\lambda_1$, $\lambda_2$, $\lambda_3$ being functions of~$r$. But the third line and (\ref{struc}) show that $\lambda_i$'s are constants and $\lambda_1=\pm\tfrac{1}{2}\lambda_2\lambda_3$. From cyclicity of $i$, $j$, $k$ we obtain two more copies of~(\ref{sdR}). Therefore,
\begin{enumerate}\itemsep=0pt
 \item[1)] $\lambda_1=\lambda_2=\lambda_3=0$ or
 \item[2)] $(\lambda_1)^2=(\lambda_2)^2=(\lambda_3)^2=4$ with $\lambda_1\lambda_2\lambda_3=\pm 8$.
\end{enumerate}
The f\/irst case leads to self-dual connection $1$-forms and Euler-top system, while the second case can be resumed to $\lambda_1=\lambda_2=\lambda_3=\pm2$ by an appropriate change of sign in~$c_i$~\cite{MR535151}. Therefore, from equations (\ref{connection1form}) and (\ref{sdR}) we get
\begin{gather*}
 \left( \frac{\partial_r c_i}{c_0} \right) = \mp \left( \frac{ c_j^2 + c_k^2 - c_i^2 }{ c_j c_k} -2\right), \qquad
 \partial_r \big( \ln c_i^2 \big) = \mp 2\big( c_j^2 + c_k^2 - c_i^2 - 2c_j c_k \big) .
 \end{gather*}
One may suppose that we must parametrize the l.h.s.\ to match the linear form in $c^2_{i}$, $c^2_{j}$ and~$c^2_{k}$ of the r.h.s.\ in the equation above. Essentially, the derivative operator aside, $c_i^2$ must be parametrized such that
\begin{gather*} \ln c_i^2= \ln \Omega_j + \ln \Omega_k - \ln \Omega_i +\text{const}= \ln \left( \frac{\Omega_j \Omega_k}{\Omega_i} \right)+\text{const}. \end{gather*}
 We choose new parametrization
 \begin{gather*}
 ( c_i )^2 = \frac{\Omega_j \Omega_k}{2\Omega_i} \ \Rightarrow \ \Omega_i = 2c_j c_k.
 \end{gather*}
 which enable us to decouple the individual parameters into their own equations turning into simpler expressions. This allows us to continue our analysis
 \begin{gather*} \partial_r \left[ \ln \left( \frac{\Omega_j \Omega_k}{2\Omega_i} \right) \right] = \frac{\dot{\Omega}_j}{\Omega_j} + \frac{\dot{\Omega}_k}{\Omega_k} - \frac{\dot{\Omega}_i}{\Omega_i} = \mp \left( \frac{\Omega_k \Omega_i}{\Omega_j} + \frac{\Omega_i \Omega_j}{\Omega_k} - \frac{\Omega_j \Omega_k}{\Omega_i} -2\Omega_i\right). \end{gather*}
 Adding up the above equation with cyclic permutations of $i$, $j$, $k$ we will f\/ind that (anti-)self-dual cases of the Bianchi~IX metric gives us
 \begin{gather*}
\frac{\dot{\Omega}_j}{\Omega_j} + \frac{\dot{\Omega}_k}{\Omega_k} - \frac{\dot{\Omega}_i}{\Omega_i} = \mp \bigg( \frac{\Omega_k \Omega_i}{\Omega_j} + \frac{\Omega_i \Omega_j}{\Omega_k} - \frac{\Omega_j \Omega_k}{\Omega_i} -2\Omega_i\bigg) \\
 \hphantom{\frac{\dot{\Omega}_j}{\Omega_j} + \frac{\dot{\Omega}_k}{\Omega_k} - \frac{\dot{\Omega}_i}{\Omega_i}}{} + \\
\frac{\dot{\Omega}_k}{\Omega_k} + \frac{\dot{\Omega}_i}{\Omega_i} - \frac{\dot{\Omega}_j}{\Omega_j} = \mp \bigg( \frac{\Omega_i \Omega_j}{\Omega_k} + \frac{\Omega_j \Omega_k}{\Omega_i} - \frac{\Omega_k \Omega_i}{\Omega_j} -2\Omega_j\bigg) \\
\hphantom{\frac{\dot{\Omega}_j}{\Omega_j} + \frac{\dot{\Omega}_k}{\Omega_k} - \frac{\dot{\Omega}_i}{\Omega_i}}{} \big{\downarrow}\\
 \Rightarrow \ \frac{\dot{\Omega}_k}{\Omega_k} = \mp 2 \left(\frac{\Omega_i \Omega_j}{\Omega_k}-\Omega_i-\Omega_j\right) \ \Rightarrow \ \dot{\Omega}_k = \mp (\Omega_i \Omega_j-\Omega_k\Omega_i-\Omega_k\Omega_j),
 \end{gather*}
where throughout derivative (denoted by dot) is taken with respect to $r$. Self-duality proceeds to give us the classical Darboux--Halphen system
\begin{gather*}
 \dot{\Omega}_i+\dot{\Omega}_j = 2\Omega_i \Omega_j .
\end{gather*}

\subsection{General Bianchi IX self-dual Einstein metric}

Following \cite{babich1998self}, we rewrite the Bianchi IX by adding a conformal scaling term $F$ in the metric
\begin{gather*}
{\rm d}s^2= F\left({\rm d}t^2+\frac{\sigma^2_1}{\Omega^2_1}+\frac{\sigma^2_2}{\Omega^2_2}+\frac{\sigma^2_3}{\Omega^2_3}\right).
\end{gather*}
where t is the cosmological time and dif\/ferent from before, here the isometry is ${\rm SU}(2)$ and ($\sigma_i$) are the corresponding ${\rm SU}(2)$ invariant forms along the spacial directions with structure constant
\begin{gather*}{\rm d}\sigma_1=\sigma_2\wedge\sigma_3,\qquad {\rm d}\sigma_2=\sigma_3\wedge\sigma_1,\qquad {\rm d}\sigma_3=\sigma_1\wedge\sigma_2.\end{gather*}
We def\/ine the new variables $A_i(t)$ by the equations
\begin{gather}
\partial_t \Omega_i=-\Omega_j\Omega_k+\Omega_i(A_j+A_k),\label{As}
\end{gather}
for distinct $i$, $j$ and $k$ taking values in the set $\{1,2,3\}$. The curvature-wise self-duality condition is expressed in terms of the new variables $A_i$ in the form of the Darboux--Halphen system
\begin{gather}\label{DarbouxA}
\partial_t A_i= -A_jA_k+A_i(A_j+A_k).
\end{gather}
Therefore we f\/ind $\Omega_i$'s by f\/irst solving system (\ref{DarbouxA}) and applying its solution in~(\ref{As}). A non-trivial solution is given by (\ref{darbouxsolution})
\begin{gather*}
A_1= 2 \frac{\partial}{\partial t}(\ln \theta_2(it)),\qquad
A_2=2 \frac{\partial}{\partial t}(\ln \theta_3(it)), \qquad
A_3=2 \frac{\partial}{\partial t}(\ln \theta_4(it)).
\end{gather*}
For simplicity, we rename $\vartheta_2\equiv\theta_2(it)$, $\vartheta_3\equiv\theta_3(it)$, $\vartheta_4\equiv\theta_4(it)$.
The system (\ref{As}) thus becomes
\begin{gather}
\partial_t{\Omega}_1=-\Omega_2\Omega_3+2\Omega_1\partial_t\ln (\vartheta_3\vartheta_4),\nonumber\\
\partial_t{\Omega}_2=-\Omega_3\Omega_1+2\Omega_2\partial_t\ln (\vartheta_4\vartheta_2),\nonumber\\
\partial_t{\Omega}_3=-\Omega_1\Omega_2+2\Omega_3\partial_t\ln (\vartheta_2\vartheta_3).\label{Omegasolution}
\end{gather}
There is a class of solutions of this system that satisf\/ies vacuum Einstein equations
\begin{gather*}
R_{ab}-\frac{1}{2}Rg_{ab}+\Lambda g_{ab}=0 ,
\end{gather*}
once we choose the appropriate conformal factor $F$ \cite{tod1994self}. This class depend on the values of the cosmological constant $\Lambda$ and satisfy the constraint
\begin{gather}\label{Omegaclass}
\vartheta^4_2\Omega_1^2-\vartheta^4_3\Omega_2^2+\vartheta^4_4\Omega_3^2=\frac{\pi^2}{4}\vartheta^4_2\vartheta^4_3\vartheta^4_4,
\end{gather}
The general two-parametric family of solutions of the system (\ref{Omegasolution}) satisfying condition (\ref{Omegaclass}), is given by the following formulas
\begin{gather*}
\Omega_1=-\frac{i}{2}\vartheta_3\vartheta_4\frac{\frac{{\rm d}}{{\rm d}q}\vartheta\big[p,q+\tfrac{1}{2}\big]}{e^{\pi i p}\vartheta[p,q]},\qquad \Omega_2=\frac{i}{2}\vartheta_2\vartheta_4\frac{\frac{{\rm d}}{{\rm d}q}\vartheta\big[p+\tfrac{1}{2},q+\tfrac{1}{2}\big]}{e^{\pi i p}\vartheta[p,q]},\nonumber\\
\Omega_3=\frac{i}{2}\vartheta_2\vartheta_4\frac{\frac{{\rm d}}{{\rm d}q}\vartheta\big[p+\tfrac{1}{2},q+\tfrac{1}{2}\big]}{e^{\pi i p}\vartheta[p,q]},
\end{gather*}
 where $\vartheta[p,q]$ denotes the theta function $\vartheta[p,q](0,ir)$, $p,q\in {\mathbb C}$. The corresponding metric is real and satisf\/ies the Einstein equations for negative cosmological constant $\Lambda$ if $p\in \mathbb{R}$ and $\mathcal{R}\{q\}=\tfrac{1}{2}$ (real part of $q$) or for positive cosmological constant if $q\in \mathbb{R}$ and $\mathcal{R}\{p\}=\tfrac{1}{2}$. In both the cases the corresponding conformal factor is given by
 \begin{gather*}
 F=\frac{2}{\pi \Lambda}\frac{\Omega_1\Omega_2\Omega_3}{\big( \frac{{\rm d}}{{\rm d}q}\ln \vartheta[p,q]\big)^2}.
 \end{gather*}
 There is another family of solutions
 \begin{gather*}
\Omega_1=\frac{1}{t+q_0}+2\frac{\partial}{\partial t}\ln\vartheta_2,\qquad \Omega_2=\frac{1}{t+q_0}+2\frac{\partial}{\partial t}\ln\vartheta_3,\qquad\Omega_3=\frac{1}{t+q_0}+2\frac{\partial}{\partial t}\ln\vartheta_4, \end{gather*}
 with $q_0\in \mathbb{R}$, that def\/ines manifolds with vanishing cosmological constant if
 \begin{gather*}F=C(t+q_0)^2\Omega_1\Omega_2\Omega_3.\end{gather*}
 \newcommand{\mat}[4]{
 \begin{pmatrix}
 #1 & #2 \\
 #3 & #4
 \end{pmatrix}
 }
\section{Gauss--Manin connection in disguise}\label{GMCD}
In this section we explain how one can derive the Darboux--Halphen equations from the Gauss--Manin connection of a universal family of elliptic curves.
This has been taken from the references
\cite{ho06-1, ho14}. The family of elliptic curves
\begin{gather*}
E_t\colon \ y^2-4(x-t_1)(x-t_2)(x-t_3)=0,\qquad t\in {\mathbb C}^3\backslash \cup_{i,j}\{t_i=t_j\},
\end{gather*}
is the universal family for the moduli of $3$-tuple $(E,(P,Q),\omega)$, where $E$ is an elliptic curve and $\omega\in H^1_{\rm dR}(E)\backslash F^1$. There is a unique regular dif\/ferential 1-form in the Hodge f\/iltration $ \omega_1\in F^1$, such that $\langle \omega,\omega_1\rangle=1$ and $\omega$, $\omega_1$ together form a basis of $H^1_{\rm dR}(E)$. $P$ and $Q$ are a pair of points of $E$ that generate the $2$-torsion subgroup with the Weil pairing $e(P,Q)=-1$. The points~$P$ and~$Q$ are given by $(t_1,0)$ and $(t_2,0)$ and $\omega=\frac{x{\rm d}x}{y}$ and $\omega_1=\frac{{\rm d}x}{y}$. The Gauss--Manin connection of the family of elliptic curves~$E_t$ written in the basis $\frac{{\rm d}x}{y}$, $\frac{x{\rm d}x}{y}$ is given as bellow
\begin{gather*}
 \nabla\begin{pmatrix}\frac{{\rm d}x}{y}\vspace{1mm}\\ \frac{x{\rm d}x}{y}\end{pmatrix}=
A
\begin{pmatrix}\frac{{\rm d}x}{y}\vspace{1mm}\\ \frac{x{\rm d}x}{y}\end{pmatrix},
\end{gather*}
where
\begin{gather*}
A =
 \frac{{\rm d}t_1}{2(t_1-t_2)(t_1-t_3)}\mat{-t_1}{1}{t_2t_3-t_1(t_2+t_3)}{t_1}\\
\hphantom{A=}{}+\frac{{\rm d}t_2}{2(t_2-t_1)(t_2-t_3)}\mat{-t_2}{1}{t_1t_3-t_2(t_1+t_3)}{t_2}\\
\hphantom{A=}{}+ \frac{{\rm d}t_3}{2(t_3-t_1)(t_3-t_2)}\mat{-t_3}{1}{t_1t_2-t_3(t_1+t_2)}{t_3}.
\end{gather*}
The reader who is not familiar with the Gauss--Manin connection must replace $\nabla$ with ${\rm d}\int_{\delta_t}$, where $t_i$'s are assumed to depend on some parameter $\tau$, ${\rm d}=\frac{\partial}{\partial \tau}$ and $\delta_t$ is a 1-dimensional homology class in~$E_t$. In the parameter space of the family of elliptic curves $E_t$ there is a unique vector f\/ield $R$, such that
\begin{gather*}
\nabla_{R}\left(\frac{{\rm d}x}{y}\right)= -\frac{x{\rm d}x}{y},\qquad \nabla_{R}\left(\frac{x{\rm d}x}{y}\right)= 0.
\end{gather*}
The vector f\/ield $R$ is given by the Darboux--Halphen system \eqref{halphen} and
it is called \textit{Gauss--Manin connection in disguise}.

\section{Frobenius manifolds and Chazy equation}\label{Frobenius}
Frobenius manifolds were developed in order to give a geometrical meaning to WDVV equations:
\begin{gather*}
 \frac{\partial^3 F(t)}{\partial t^{\alpha} \partial t^{\beta} \partial t^{\lambda}} \eta^{\lambda \mu} \frac{\partial^3 F(t)}{\partial t^{\mu} \partial t^{\gamma} \partial t^{\delta}} =
 \frac{\partial^3 F(t)}{\partial t^{\delta} \partial t^{\beta} \partial t^{\lambda}} \eta^{\lambda \mu} \frac{\partial^3 F(t)}{\partial t^{\mu} \partial t^{\gamma} \partial t^{\alpha}},
\end{gather*}
where $F(t)$, with $t=(t^1,t^2,\dots,t^n)$, is a quasi-homogeneous function on its parameters.
The above equations conceal properties of an associative commutative algebra on the tangent space of a~manifold $M$ of dimension $n$
def\/ined by the parameter space $(t^1,t^2,\dots,t^n)$. That's the essence of a Frobenius manifold that we will
detail below starting with the algebraic structure in $TM$.

\subsection{Frobenius algebra}

An algebra $A$ over $\mathbb{C}$ is Frobenius if
\begin{itemize}
\item it is a commutative associative $\mathbb{C}$-algebra with unity $e$,
\item it has a $\mathbb{C}$-bilinear symmetric non-degenerate inner product
\begin{align*}
\langle \, , \, \rangle\colon \ & A\times A\longrightarrow \mathbb{C}, \\
& (a , b)\mapsto \langle a,b \rangle,
\end{align*}
which is invariant, i.e., $\langle a.b,c \rangle =\langle a,b.c \rangle$
\end{itemize}
{\bf Properties:}
Let $e_\alpha$, $\alpha=1,\dots, N$, be any basis in $A$, such that $e_1=e$ is the unity. By notation, we def\/ine $\eta_{\alpha\beta}:=\langle e_\alpha,e_\beta\rangle$, which yields the matrix $\eta:=[\eta_{\alpha\beta}]_{1\leq \alpha, \beta \leq N}$ and its inverse $\eta^{-1}:=[\eta^{\alpha\beta}]_{1\leq \alpha, \beta \leq N}$, and it follows $\eta^{\alpha\beta}\eta_{\beta\gamma}=\delta^{\alpha}_{\gamma}$. By writing $e_\alpha\cdot e_\beta$ in the given basis, we f\/ind the structure constants $c_{\alpha\beta}^\gamma$ def\/ined by $e_\alpha \cdot e_\beta= c_{\alpha\beta}^\gamma e_\gamma$. If we set $c_{\alpha\beta\gamma}=c_{\alpha\beta}^\epsilon\eta_{\epsilon\gamma}$, then we get $c_{\alpha\beta}^\gamma=c_{\alpha\beta\epsilon}\eta^{\epsilon\gamma}$.
Note that in all above expressions, and in what follows, Einstein summation of indices is implicit.
Therefore, $\eta_{\alpha\beta}$ and the structure constants $c_{\alpha\beta}^\gamma$ satisfy
\begin{alignat}{3}\label{commutativity}
& \textnormal{commutativity} \qquad&& \eta_{\alpha\beta}=\eta_{\beta\alpha},& \\
& \textnormal{associativity}\qquad&& (e_\alpha. e_\beta).e_\gamma=e_\alpha.( e_\beta .e_\gamma) \therefore c_{\alpha\beta}^\epsilon c_{\epsilon\gamma}^\delta=c_{\alpha\epsilon}^\delta c_{\delta\beta}^\epsilon,& \label{associativity}\\
& \textnormal{normalization}\qquad&& c_{1\beta}^\alpha=\delta_\beta^\alpha,& \\ \label{invariance}
& \textnormal{invariance\&commutat.}\qquad&& c_{\alpha\beta\gamma}= \langle e_\alpha e_\beta, e_\gamma\rangle=c_{\beta\alpha\gamma}=c_{\alpha\gamma\beta}.&
\end{alignat}
Now consider an $n$-parametric deformation of the Frobenius algebra $A_t$, $t= (t^1,t^2,\dots,t^n)$, with structure constants $c_{\alpha\beta}^\gamma(t)$ preserving relations \eqref{commutativity} to \eqref{invariance}. Such deformed algebra~$A_t$ can be seen as a f\/iber bundle with the space of parameters $t\in M$ as base space. We identify this f\/iber bundle with the tangent bundle $TM$ to arrive at the def\/inition of a Frobenius manifold. The requirements for this to happen are presented in the def\/inition below.

\subsection{Frobenius manifold}

A Frobenius manifold M of dimension $n$, is an $n$-dimensional Riemannian manifold, such that for all $t\in M$ the tangent space $T_tM$ contains the structure of a Frobenius algebra $ (A_t, \langle \,, \,\rangle_t)$,
satisfying the following axioms:
\begin{enumerate}\itemsep=0pt
 \item[A.1.] The metric $\langle \,,\, \rangle_t$ on $M$ is f\/lat. The unit vector $e$ must be f\/lat, i.e., $\nabla e = 0$, where $\nabla$ is the Levi-Civita connection for the
metric.
 \item[A.2.] Let $c$ be the 3-tensor $c(x,y,z) = \langle x.y ,z \rangle $, with $x, y, z \in T_tM$. Then the 4-tensor $(\nabla_w c)(x,$ $y,z)$ must be symmetric in $x,y,z,w \in T_tM$.
 \item[A.3.] A linear vector f\/ield $E$ must be f\/ixed on $M$, i.e., $\nabla( \nabla E) = 0$ such that the corresponding one-parameter group of dif\/feomorphisms acts by conformal
transformations of the metric $\langle \, , \, \rangle$ and by rescaling on the Frobenius algebras $T_tM$.
\end{enumerate}
The f\/latness of the metric $\langle \,, \,\rangle$ implies the existence of a system of f\/lat coordinates $t^1, \dots, t^n$ on~$M$. In these f\/lat coordinates the structure
constants of $A_t$ are given by
\begin{gather*}
\frac{\partial}{\partial t^{\alpha}}.\frac{\partial}{\partial t^{\beta}} = c_{\alpha \beta}^{\gamma}(t) \frac{\partial}{\partial t^{\gamma}}.
\end{gather*}

{\bf Potential deformation.} If there is a function $F(t)$, called {\it potential}, such that the structure constants of $A_t$, $t\in M$, can be locally represented as
\begin{gather*}
 c_{\alpha \beta\gamma}(t) = \frac{\partial^3 F(t)}{\partial t^{\gamma} \partial t^{\alpha} \partial t^{\beta}},
\end{gather*}
 satisfying A.2 with unity vector $e=\tfrac{\partial}{\partial t^1}$, and the metric given by
\begin{gather*}
 \eta_{\beta\gamma}=c_{1 \beta\gamma} = \frac{\partial^3 F(t)}{\partial t^{\gamma} \partial t^{1} \partial t^{\beta}}\qquad \textnormal{s.t.}\qquad \frac{\partial^4 F(t)}{\partial t^{\alpha}\partial t^{\gamma} \partial t^{1} \partial t^{\beta}}=0,
\end{gather*}
 satisfying A.1, and the associativity property (\ref{associativity}) represented by the WDVV equations
\begin{gather*}
 \frac{\partial^3 F(t)}{\partial t^{\alpha} \partial t^{\beta} \partial t^{\lambda}} \eta^{\lambda \mu} \frac{\partial^3 F(t)}{\partial t^{\mu} \partial t^{\gamma} \partial t^{\delta}} =
 \frac{\partial^3 F(t)}{\partial t^{\delta} \partial t^{\beta} \partial t^{\lambda}} \eta^{\lambda \mu} \frac{\partial^3 F(t)}{\partial t^{\mu} \partial t^{\gamma} \partial t^{\alpha}},
\end{gather*}
 then $M$ is a Frobenius manifold
 and the Frobenius algebra $A_t$ is called a potential deformation.
 Note that the condition~A.3 is satisf\/ied by a quasihomogeneous function $F(t)$.
\begin{Example}

Let $\dim M = 3$, and consider the basis $e = e_1=\tfrac{\partial}{\partial t^1}$, $e_2=\tfrac{\partial}{\partial t^2}$ and $e_3=\tfrac{\partial}{\partial t^3}$ of the 3-dimensional algebra~$A_t$. Then the multiplication law is given by
\begin{gather*}
e_2^2 = f_{xxy}e_1 + f_{xxx}e_2 + e_3, \qquad
e_2 e_3 = f_{xyy}e_1 + f_{xxy}e_2, \qquad
e_3^2 = f_{yyy}e_1 + f_{xyy}e_2,
\end{gather*}
where the funtion $F(t)$ has the form $F(t) = \frac{1}{2} (t^1)^2t^3 + \frac{1}{2} t^1(t^2)^2 + f(t^2,t^3)$ and the notation $f_x= \partial_x f(x,y)$, $f_y= \partial_y f(x,y)$. The associativity condition $(e_2^2) e_3 = e_2 (e_2 e_3)$
implies the following PDE for $f(x,y)$:
\begin{gather}\label{associativity3dim}
 f^2_{xxy} = f_{yyy} + f_{xxx}f_{xyy}.
\end{gather}
\end{Example}

\subsection{Chazy equation and Darboux--Halphen system}\label{modular frobenius}
In this section we explain how the Chazy equation arises from a 3-dimensional Frobenius manifold. We follow Dubrovin's notes \cite{Dubrovin1992NuPhB,Dubrovin:1996di}.
Let dim $M = 3$, and consider the potential function
\begin{gather*}
F(t)= \frac{1}{2} \big(t^1\big)^2t^3 + \frac{1}{2} t^1\big(t^2\big)^2 - \frac{\big(t^2\big)^4}{16}\gamma\big(t^3\big),
\end{gather*}
where $\gamma(\tau)$ is an unknown $2\pi$-periodic function that is analytic at $\tau=i \infty$. Then the associativity condition~\eqref{associativity3dim} leads to the Chazy equation
\begin{gather*}
\gamma'''= 6\gamma\gamma''-9(\gamma')^2.
\end{gather*}
The solution, up to a shift in $\tau$, is given by $\gamma(\tau)= \frac{\pi i}{3}E_2(\tau)$, where $E_2$ is the weight-2 Eisenstein series.
Notice that the Darboux--Halphen solution (\ref{darbouxsolution}) leads to
\begin{gather*}
t_1+t_2+t_3= \frac{\pi i}{2}E_2(\tau),
\end{gather*}
which can be easily checked from \eqref{27spet2017} or by
writing the theta functions in terms of Dedekind eta
function, see \cite[Chapter~3, p.~29]{Zagier:2008uf}.
Applying $\tau$ derivatives on both sides and using Darboux--Halphen equations, one can also check that the solution to Darboux--Halphen system~(\ref{darbouxsolution}) are the roots of the cubic equation
\begin{gather*}
y^3-\frac{3}{2}\gamma(\tau)y^2+\frac{3}{2}\gamma'(\tau)y -\frac{1}{4}\gamma''(\tau)=0.
\end{gather*}

\section{Conclusion}
The study of Darboux--Halphen equations in several dif\/ferent problems in theoretical physics and mathematics raised more and
more questions that eventually lead us to further studies.

The problem involving Gauss--Manin connection in disguise lies at the center of some questions. It shows that the
Darboux--Halphen system corresponds to a vector f\/ield in the moduli of an enhanced elliptic curve.
As mentioned in Section~\ref{BianchiIXsection}, the Bianchi~IX four-manifold~\eqref{bianchiix} also describes the
reduced moduli of 2-monopoles and its self-dual curvature equations can be reparametrized to the Darboux--Halphen equations.
Furthermore, in the problem of 2-monopoles it has been found that a 2-monopole solution relates to an elliptic curve as its
spectral curve \cite{1982CMaPh..83..579H,Hitchin:1983ay,Hurtubise:1983hl}. Therefore, we believe that Gauss--Manin connection in disguise is a new way to demonstrate the association of spectral curves and the curvature equations of the moduli of monopole solutions. Starting from these coincidences, in~\cite{Torres:2017hqc} one of the authors started to f\/ind more evidences to support this idea.

Another interesting remark is the fact that potential functions and structure constants in Frobenius manifolds correspond
to prepotentials (or genus zero topological partition function) and Yukawa couplings in topological string theory and
Gauss--Manin connection in disguise has been used in the moduli of enhanced Calabi--Yau varieties to f\/ind polynomial expressions for
Yukawa couplings and higher genus topological partition functions \cite{Alim:2014dea, 2011arXiv1111.0357M}.
It would be interesting to f\/ind cases where the moduli of enhanced Calabi--Yau varieties are also Frobenius manifolds.
In particular, the Frobenius manifold presented in Section~\ref{modular frobenius} is a case of modular Frobenius
manifold where the prepotential is preserved under a inverse symmetry that acts as an $S$ gene\-ra\-tor of the modular
group ${\rm SL}(2,\mathbb{Z})$ in $t^3$ direction \cite{Morrison:2011cd}. Such modularity is a~desirable property
that can establish a relation to Gauss--Manin connection in disguise and may be
extended to
the group of transformations of Calabi--Yau modular forms \cite{Alim:2014dea,2011arXiv1111.0357M}.

\subsection*{Acknowledgements}

During the period of preparation of the manuscript MACT was fully sponsored by CNpQ-Brasil. The research of RR was supported by FAPESP through Instituto de Fisica, Universidade de Sao Paulo with grant number 2013/17765-0. The work was initiated during the visit of RR to IMPA, he would like to thank IMPA for the hospitality during the course of this project.

\pdfbookmark[1]{References}{ref}
\LastPageEnding

\end{document}